\documentclass[10pt]{report}
\usepackage{amssymb}
\usepackage{latexsym}
\usepackage[dvips]{graphicx}
\usepackage{natbib}
\usepackage{amsmath}
\usepackage{float}
\usepackage{labelfig}
\usepackage{epsfig,psfig}

\def\lemma{\textbf{Lemma}}
\def\theorem{\textbf{Theorem}}

\def\corollary{\textbf{Corollary}}

\def\diamond{\diamondsuit}

\newcommand{\Map}{\operatorname{Map}}

\begin{document}

\begin{center}
\textbf{Constructing convex planes in the pants complex.}

Javier Aramayona, Hugo Parlier, Kenneth J. Shackleton

[First draft: February 2007]\\~\\

\end{center}

\setlength{\parindent}{0em}

A\begin{small}BSTRACT\end{small}: Our main theorem identifies a class of totally geodesic subgraphs of the $1$-skeleton of the pants complex, each isomorphic to the product of two Farey graphs. We deduce the existence of many convex planes in the $1$-skeleton of the pants complex.\\

K\begin{small}EYWORDS\end{small}: pants complex; Weil-Petersson metric\\



\setlength{\parindent}{0em}

\textbf{$\S1$. Introduction.\\}

Let $\Sigma$ be a compact, connected and orientable surface, possibly with non-empty boundary, of genus $g(\Sigma)$ and $|\partial \Sigma|$ boundary components, and refer to as the \textit{mapping class group} $\Map(\Sigma)$ the group of all self-homeomorphisms of $\Sigma$ up to homotopy.

\setlength{\parindent}{2em}

After Hatcher-Thurston [HT], to the surface $\Sigma$ one may associate a simplicial graph $\mathcal{P}(\Sigma)$, the \textit{pants graph}, whose vertices are all the pants decompositions of $\Sigma$ and any two vertices are connected by an edge if and only if they differ by an elementary move; see $\S2.2$ for an expanded definition. This graph is connected, and one may define a path-metric $d$ on $\mathcal{P}(\Sigma)$ by first assigning length $1$ to each edge and then regarding the result as a length space.

The pants graph, with its own geometry, is a fundamental object to study. Brock [B] revealed deep connections with hyperbolic $3$-manifolds and proved the pants graph is the correct combinatorial model for the Weil-Petersson metric on Teichm\"uller space, for the two are quasi-isometric. The isometry group of $(\mathcal{P}, d)$ is also correct in so far as the study of surface groups is concerned, for Margalit [Mar] proved it is almost always isomorphic to the mapping class group of $\Sigma$. In addition, Masur-Schleimer [MasS] proved the pants graph of any closed surface of genus at least $3$ to be one-ended. With only a few exceptions, the pants graph is not hyperbolic in the sense of Gromov [BF].

In [APS], the authors prove that every subgraph of $\mathcal{P}$ isomorphic to the Farey graph is in fact totally geodesic in $(\mathcal{P}, d)$. The purpose of this paper is to study the extrinsic geometry of another class of subgraphs of the pants graph, each determined by $2$-handle multicurves as defined at the end of $\S2.1$.\\

\newtheorem{1}{\theorem}

\begin{1}
Let $\Sigma$ be a compact, connected and orientable surface, and denote by $Q$ any $2$-handle multicurve on $\Sigma$. Then, $\mathcal{P}_{Q}$ is totally geodesic in $\mathcal{P}(\Sigma)$.\\
\label{main}
\end{1}

The completion of the Weil-Petersson metric can be characterised by attaching so-called strata [Mas]. These are totally geodesic subspaces of the completion, by a result of Wolpert [W], and correspond to lower dimensional Teichm\"uller spaces, or products thereof, each with their own Weil-Petersson metric, or product metric. Combining this with Theorem 1.1 of Brock [B], one finds the $2$-handle subgraphs of the pants graph are uniformly quasi-convex.  Still, Theorem \ref{main} is not implied by any known coarse geometric result. Moreover, Theorem \ref{main} establishes a complete analogy between the geometry of the $2$-handle subgraphs in a pants graph and the geometry of the corresponding strata lying in the completed Weil-Petersson space.

In order to prove Theorem \ref{main}, we shall project paths in the pants graph to paths in the given $2$-handle graph of no greater length. All the notation of Theorem \ref{project} will be explained in $\S2$, but for now we point out the finite set of curves $\pi_{Q}(\nu)$ is the subsurface projection after Masur-Minsky [MasMi] of a pants decomposition $\nu$ to the subgraph $\mathcal{P}_{Q}$ determined by the $2$-handle multicurve $Q$. Note, our definition differs slightly from that given in [APS]. The intrinsic metric on the graph $\mathcal{P}_{Q}$, assigning length $1$ to each edge, is denoted by $d_{Q}$.\\

\newtheorem{2}[1]{\theorem}

\begin{2}
Let $\Sigma$ be a compact, connected and orientable surface and denote by $Q$ any $2$-handle multicurve on $\Sigma$. Let $(\nu_{0}, \ldots, \nu_{n})$ be a path in the pants graph $\mathcal{P}(\Sigma)$ such that $\nu_{n}$ is a vertex of $\mathcal{P}_{Q}$. For each index $i \leq n - 1$ and for each $\omega_{i} \in \pi_{Q}(\nu_{i})$, there exists an integer $j \in \{1, 2\}$ and a pants decomposition $\omega_{i+j} \in \pi_{Q}(\nu_{i+j})$ of $\Sigma$ such that $d_{Q}(\omega_{i}, \omega_{i+j}) \leq j$.\\
\label{project}
\end{2}

To the authors' knowledge, it has yet to be decided whether there exists a distance non-increasing projection from the whole pants graph to any one of its $2$-handle subgraphs. In the absence of an affirmative result, Theorem \ref{project} may well hold independent interest.

Let us indicate two consequences of Theorem \ref{main}. First, by considering a pair of bi-infinite geodesics, one in either factor Farey graph for a $2$-handle multicurve, we deduce the following. By a \textit{plane} we shall mean a graph isomorphic to the Cayley graph of the group $\Bbb{Z} \oplus \Bbb{Z}$ with standard generating set.\\

\newtheorem{3}[1]{\corollary}

\begin{3}
Let $\Sigma$ be a compact, connected and orientable surface of complexity at least $3$. Then, $\mathcal{P}(\Sigma)$ contains infinitely many convex planes.\\
\label{planes}
\end{3}

Second, we exhibit convex planes in the pants graph invariant under the action of a particular family of mapping classes.\\

\newtheorem{3a}[1]{\corollary}

\begin{3a}
Let $f \in \Map(\Sigma)$ be a mapping class fixing two disjoint and incompressible complexity $1$ subsurfaces of $\Sigma$, acting on each as a pseudo-Anosov mapping class and fixing a pants decomposition of their complement. Then, there exists a convex plane in $\mathcal{P}(\Sigma)$ on which $f^{2}$ acts by translation.\\
\label{invplane}
\end{3a}

The plan of this paper is as follows. In $\S2$ we recall all the terminology we need, much of which is already standard. In $\S3$ we give an elementary proof to Theorem \ref{project}. In $\S4$ we apply Theorem \ref{project} to give an elementary proof to Theorem \ref{main}. Finally, in $\S5$ we prove Corollary \ref{planes}.\\

\setlength{\parindent}{0em}

A\begin{small}CKNOWLEDGEMENTS:\end{small} The first author was supported by a short term research fellowship at the Universit\'e de Provence, and is grateful for both its financial support and warm hospitality. The third author was supported by a long term Japan Society for the Promotion of Science postdoctoral fellowship, number P06034, and he gratefully acknowledges the JSPS for its financial support. The third author gratefully acknowledges the Department of Mathematical and Computing Sciences at the Tokyo Institute of Technology for its warm hospitality. The authors wish to thank Saul Schleimer for an interesting conversation.\\~\\

\textbf{$\S2$. Background and definitions.}\\

We supply all the background and terminology needed both to understand the statements of our main results, and to make sense of their proofs. Throughout, we define a \textit{loop} on $\Sigma$ as the homeomorphic image of a standard circle. A subsurface of $\Sigma$ is said to be \textit{incompressible} only if its inclusion descends to an injection on fundamental groups.\\

\setlength{\parindent}{0em}

\textbf{$\S2.1.$ Curves and multicurves.} A loop on $\Sigma$ is said to be \textit{trivial} only if it bounds a disc and \textit{peripheral} only if it bounds an annulus whose other boundary component belongs to $\partial \Sigma$. For a non-trivial and non-peripheral loop $c$, we shall denote by $[c]$ its free homotopy class. A \textit{curve} is by definition the free homotopy class of a non-trivial and non-peripheral loop. Given any two curves $\alpha$ and $\beta$, their intersection number $\iota(\alpha, \beta)$ is defined equal to $\rm{min}\{|a \cap b| : a \in \alpha, b \in \beta\}$.

\setlength{\parindent}{2em}

We shall say two curves are \textit{disjoint} only if they have zero intersection number, and otherwise say they \textit{intersect essentially}. A pair of curves $\{\alpha, \beta\}$ is said to \textit{fill} the surface $\Sigma$ only if $\iota(\delta, \alpha) + \iota(\delta, \beta) > 0$ for every curve $\delta$. In other words, every curve on $\Sigma$ intersects at least one of $\alpha$ and $\beta$ essentially.

A \textit{multicurve} is a collection of distinct and disjoint curves, and the intersection number for a pair of multicurves is to be defined additively. We denote by $\kappa(\Sigma)$ the cardinality of any maximal multicurve on $\Sigma$, equal to $3g(\Sigma) + |\partial \Sigma| - 3$, and refer to this as the \textit{complexity of $\Sigma$}. Note, the only surfaces of complexity $1$ are the $4$-holed sphere and the $1$-holed torus.

Given a set of disjoint loops $L$, such as the boundary of some subsurface of $\Sigma$, we denote by $[L]$ the multicurve maximal among all multicurves whose every curve is represented by some element of $L$. We shall say a multicurve $\omega$ has \textit{codimension $k$}, for some non-negative integer $k$, only if $|\omega| = \kappa(\Sigma) - k$. We shall say a codimension $2$ multicurve $Q$ is a \textit{$2$-handle multicurve} only if the complement of every simple representative of $Q$ contains two complexity $1$ components, each either a $1$-holed torus or a $4$-holed sphere containing three components of $\partial \Sigma$.\\


\setlength{\parindent}{0em}

\textbf{$\S2.2.$ Pants decompositions.} A \textit{pants decomposition} of a surface is a maximal collection of distinct and disjoint curves, in other words a maximal multicurve. Two pants decompositions $\mu$ and $\nu$ are said to be related by an \textit{elementary move} only if $\mu \cap \nu$ is a codimension $1$ multicurve and the remaining two curves together either fill a $4$-holed sphere and intersect twice or fill a $1$-holed torus and intersect once; consider Figure \ref{fig:elementarymoves} below.\\~\\

\begin{figure}[h]
\begin{center}
\AffixLabels{\centerline{\epsfig{file = 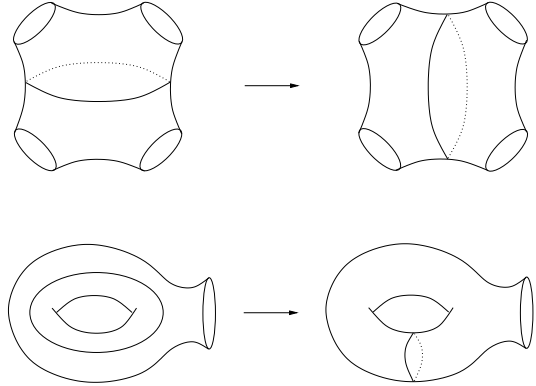, width=7cm, angle= 0}}}
\end{center}
\caption{The two types of elementary move.}
\label{fig:elementarymoves}
\end{figure}

\textbf{$\S2.3.$ Arcs.} An \textit{arc} on $\Sigma$ is the homotopy class, relative to $\partial \Sigma$, of an embedded interval ending on $\partial \Sigma$ that does not bound a disc with $\partial \Sigma$. There are broadly two types of arc: those that end on only one component of $\partial \Sigma$, referred to as \textit{waves}, and those that end on two different components of $\partial \Sigma$, referred to as seams. In this paper, we shall be considering only waves.

\setlength{\parindent}{2em}

We may similarly define the intersection number of a pair of arcs, or an arc and a curve, and say two arcs are disjoint or intersect essentially.\\

\setlength{\parindent}{0em}

\textbf{$\S2.4.$ Graphs and paths.} For us, a \textit{path} in a graph shall be a finite sequence of vertices such that any consecutive pair spans an edge; one can recover a topological path by joining up the dots. A \textit{geodesic} is then a path realising distance. A subgraph $F$ of a graph $G$ is said to be \textit{totally geodesic} only if every geodesic in $G$ whose two endpoints belong to $F$ is entirely contained in $F$. A subgraph $F$ of a graph $G$ is said to be \textit{convex} only if any two vertices of $F$ are connected by a geodesic entirely contained in $F$. Any totally geodesic subgraph is convex, though a convex subgraph need not be totally geodesic. Finally, the \textit{product of two graphs $G_{1}$ and $G_{2}$}, denoted $G_{1} * G_{2}$, is the graph with vertex set $G_{1} \times G_{2}$ and two distinct vertices $(x_{1}, x_{2})$ and $(y_{1}, y_{2})$ are connected by an edge if and only if either $x_{1}$ and $y_{1}$ span an edge in $G_{1}$ or $x_{2}$ and $y_{2}$ span an edge in $G_{2}$.\\

\textbf{$\S2.5.$ Subsurface projections.} Given a curve $\alpha$ and an incompressible subsurface $Y$ of $\Sigma$, we shall write $\alpha \subset Y$ only if $\alpha$ can be represented by a non-peripheral loop on $Y$. If every loop representing $\alpha$ has non-empty intersection with $Y$ we can say \textit{$\alpha$ and $Y$ intersect}, otherwise we say they are \textit{disjoint}. If every loop representing $\alpha$ intersects $Y$ in at least one interval, we can say \textit{$\alpha$ crosses $Y$}.

\setlength{\parindent}{2em}

Let $Y$ denote any complexity $1$ incompressible subsurface of $\Sigma$. Let $\alpha$ be any curve intersecting $Y$, and choose any simple representative $c \in \alpha$ such that $\#(c \cap \partial Y)$ is minimal. We refer to each component of $c \cap Y$ as a \textit{footprint of $c$ on $Y$}, and to the homotopy class, relative to $\partial Y$, of such a footprint as a \textit{footprint of $\alpha$ on $Y$}. Note, footprints of a curve can be arcs or curves.

Given a footprint $b$ for the curve $\alpha$ there only ever exists one curve on $Y$ disjoint from $b$, and such a curve shall be referred to as a \textit{projection of $\alpha$}. Note the set of $\alpha$ projections, each counted once, depends only on $\alpha$ and the isotopy class of the surface $Y$, and we denote this set by $\pi_{Y}(\alpha)$. For a second multicurve $\nu$ we may similarly define $\pi_{Y}(\nu)$. If $\nu$ is disjoint from $Y$ then we define $\pi_{Y}(\nu)$ to the empty set. The set $\pi_{Y}(\nu)$ is well-defined, and is an example of a \textit{subsurface projection} as defined by Masur-Minsky in $\S1.1$ of [MasMi].

By way of example, we note that $\pi_{Y}(\delta) = \{\delta\}$ for any curve $\delta \subset Y$. If $\delta \subset Y$ is a curve and $\alpha$ is a second curve crossing $Y$ and disjoint from $\delta$, then $\delta \in \pi_{Y}(\alpha)$.

Finally, if $Q$ is a $2$-handle multicurve and $Y_{1}$ and $Y_{2}$ are two non-isotopic incompressible complexity $1$ subsurfaces of $\Sigma$ disjoint from $Q$, for any multicurve $\nu$ on $\Sigma$ we define $\pi_{Q}(\nu)$ to be equal to $\{\{\delta_{1}, \delta_{2}\} \cup Q : \delta_{1} \in \pi_{Y_{1}}(\nu), \delta_{2} \in \pi_{Y_{2}}(\nu)\}$. It follows $\pi_{Q}(\nu)$ is the empty set whenever $\nu$ is disjoint from both $Y_{1}$ and $Y_{2}$. However, if $\nu$ is also a pants decomposition, then each element of $\pi_{Q}(\nu)$ is a pants decomposition containing $Q$ and thus is a vertex of $\mathcal{P}_{Q}$. Moreover, $\pi_{Q}$ restricts to the identity on $\mathcal{P}_{Q}$.\\~\\

\setlength{\parindent}{0em}

\textbf{$\S3.$ One proof of Theorem 2.}\\

\newtheorem{4}[1]{\lemma}

\setlength{\parindent}{2em}

\noindent We begin with three results, the third of which plays an especially important role in the proof of Theorem \ref{project}. In what follows, we shall make use of the fact that a pair of disjoint waves projects to a pair of curves either equal or intersecting minimally.

\begin{4}
Let $Q$ be a codimension $1$ multicurve on $\Sigma$, and denote by $Y$ an incompressible complexity $1$ subsurface of $\Sigma$ disjoint from $Q$. For two disjoint waves or curves $a_{1}$ and $a_{2}$ on $Y$, denote by $\alpha_{1}$ and $\alpha_{2}$ the unique curves on $Y$ such that $\iota(\alpha_{1}, a_{1}) = 0$ and such that $\iota(\alpha_{2}, a_{2}) = 0$. Then, $d(\{\alpha_{1}\} \cup Q, \{\alpha_{2}\} \cup Q) \leq 1$.
\label{smalldist}
\end{4}

\noindent \textbf{Proof:} If at least one of $a_{1}$ and $a_{2}$ is a curve, then $\alpha_{1}$ and $\alpha_{2}$ are equal and $d(\{\alpha_{1}\} \cup Q, \{\alpha_{2}\} \cup Q) = 0$. Otherwise, $\alpha_{1}$ and $\alpha_{2}$ are either equal or intersect minimally, and as such $d(\{\alpha_{1}\} \cup Q, \{\alpha_{2}\} \cup Q) \leq 1$. $\diamond$\\

\newtheorem{5}[1]{\lemma}

\begin{5}
Let $Q$ be a $2$-handle multicurve on $\Sigma$, and let $\nu_{0}$ and $\nu_{1}$ be two vertices of $\mathcal{P}(\Sigma)$ such that $d(\nu_{0}, \nu_{1}) = 1$. For $\omega_{0} \in \pi_{Q}(\nu_{0})$, if there exists $\omega_{1} \in \pi_{Q}(\nu_{1})$ such that $\omega_{0} \cap \omega_{1}$ is not equal to $Q$, then there exists $\omega_{1}' \in \pi_{Q}(\nu_{1})$ such that $d(\omega_{0}, \omega_{1}') \leq 1$.
\label{goodproj}
\end{5}

\noindent \textbf{Proof:} Let $\omega_{0} \in \pi_{Q}(\nu_{0})$, and supppose there exists $\omega_{1} \in \pi_{Q}(\nu_{1})$ such that $\omega_{0} \cap \omega_{1}$ is not equal to $Q$. Note then, $\omega_{0} \cap \omega_{1}$ is a codimension $1$ multicurve properly containing $Q$. If $d(\omega_{0}, \omega_{1}) \geq 2$, then necessarily there exists a pants decomposition $\omega_{1}' \in \pi_{Q}(\nu_{0} \cap \nu_{1})$. Since $\omega_{1}' \in \pi_{Q}(\nu_{0})$, so $d(\omega_{0}, \omega_{1}') \leq 1$. If on the other hand $d(\omega_{0},\omega_{1}) \leq 1$, then we may define $\omega_{1}'$ to be equal to $\omega_{1}$. $\diamond$\\

\newtheorem{6}[1]{\lemma}

\begin{6}
Let $Q$ be a $2$-handle multicurve on $\Sigma$, and denote by $Y_{1}$ and $Y_{2}$ two non-isotopic incompressible complexity $1$ subsurfaces of $\Sigma$ disjoint from $Q$. Let $\nu_{0}, \nu_{1}, \nu_{2}$ be any geodesic in $\mathcal{P}(\Sigma)$ of length $2$, and let $\omega_{0} \in \pi_{Q}(\nu_{0})$. If $\omega_{0} \cap \omega_{1}$ is equal to $Q$ for each $\omega_{1} \in \pi_{Q}(\nu_{1})$, then the multicurve $\nu_{0} \cap \nu_{1} \cap \nu_{2}$ intersects both $Y_{1}$ and $Y_{2}$.
\label{trick}
\end{6}

\noindent \textbf{Proof:} There exists a unique curve $\delta \in \nu_{0}$ such that $\omega_{0} \in \pi_{Q}(\delta)$. Let $R$ denote the multicurve $\nu_{0} \cap \nu_{1} \cap \nu_{2}$, noting $|R| \geq \kappa(\Sigma) - 2$.

Suppose for contradiction that $R$ does not intersect $Y_{1}$. Then, we may project $\delta$ to the complement of $Y_{1}$ in $\Sigma$ to find a curve $\delta'$ disjoint from $Y_{1}$ such that $\delta'$ intersects $Y_{2}$ and such that $\iota(\omega_{0} - Q, \delta')$ is zero.

We now note that $\delta'$ is distinct and disjoint from every curve in $R$, for respectively $\omega_{0} \cap \omega$ is equal to $Q$ for each $\omega \in \pi_{Q}(R)$, by assumption, and both $\iota(R, [\partial Y_{1}])$ and $\iota(R, \delta)$ are zero. It follows that $R \sqcup \{\delta'\}$ is a multicurve disjoint from $Y_{1}$ and, since it cannot contain $[\partial Y_{1}]$, as such has cardinality at most $\kappa(\Sigma) - 2$. Thus, $$|R| = |(R \sqcup \{\delta'\}) - \{\delta'\}| = |R \sqcup \{\delta'\}| - |\{\delta'\}| \leq \kappa(\Sigma) - 2 - 1 = \kappa(\Sigma) - 3.$$ To be more succinct, $|R| \leq \kappa(\Sigma) - 3$. We therefore have two incompatible estimates for the cardinality of $R$, and this is a contradiction.

A parallel argument applies to $Y_{2}$, and the statement of the lemma therefore holds. $\diamond$\\

\setlength{\parindent}{2em}

We now turn to proving Theorem \ref{project}, denoting by $Y_{1}$ and $Y_{2}$ two incompressible non-isotopic complexity $1$ subsurfaces of $\Sigma$ disjoint from $Q$. Let $\omega_{0}$ denote any element of $\pi_{Q}(\nu_{0})$.

Suppose inductively we have chosen the vertex $\omega_{k} \in \pi_{Q}(\nu_{k})$, for some $k \geq 0$. If there exists a pants decomposition $\omega \in \pi_{Q}(\nu_{k+1})$ such that $\omega_{k} \cap \omega$ is not equal to $Q$, then by Lemma \ref{goodproj} there exists a pants decomposition $\omega' \in \pi_{Q}(\nu_{k+1})$ such that $d(\omega_{k}, \omega') \leq 1$. We define $\omega_{k+1}$ to be equal to $\omega'$.

We now consider the remaining case, that $\omega_{k} \cap \omega$ is equal to $Q$ for every pants decomposition $\omega \in \pi_{Q}(\nu_{k+1})$. Note then, $k \leq n-2$. By Lemma \ref{trick}, there exists a multicurve $R$, contained in $\nu_{k} \cap \nu_{k+1} \cap \nu_{k+2}$, such that $R$ intersects both $Y_{1}$ and $Y_{2}$. Taking one footprint of $R$ on $Y_{1}$ and then on $Y_{2}$, we may construct vertices $\omega_{k+1}$ and $\omega_{k+2}$ of $\mathcal{P}_{Q}$ such that $\omega_{k} \cap \omega_{k+1}$ and $\omega_{k+1} \cap \omega_{k+2}$ are both codimension $1$ multicurves, and where $\omega_{k+2} \in \pi_{Q}(R) \subseteq \pi_{Q}(\nu_{k+2})$. Note though, $\omega_{k+1}$ need not be contained in $\pi_{Q}(\nu_{k+1})$. As each footprint of $\omega$ on $Y_{1}$ and on $Y_{2}$ is either a wave or a curve, by Lemma \ref{smalldist} we have in turn $d(\omega_{k}, \omega_{k+1}) \leq 1$ and $d(\omega_{k+1}, \omega_{k+2}) \leq 1$. Thus, $d(\omega_{k}, \omega_{k+2}) \leq 2$ and the induction continues from $k+2$.

This concludes a proof of Theorem \ref{project}. $\diamond$\\\\

\setlength{\parindent}{0em}

\textbf{$\S4$. One proof of Theorem \ref{main}.}\\

Let $Q$ be a $2$-handle multicurve on $\Sigma$. Suppose, for contradiction, that $\mathcal{P}_{Q}$ is not totally geodesic. Then, there exist two vertices $\mu$ and $\nu$ of $\mathcal{P}_{Q}$ and a geodesic $\mu = \nu_{0}, \nu_{1}, \ldots, \nu_{n} = \nu$ in $\mathcal{P}(\Sigma)$ not entirely contained in $\mathcal{P}_{Q}$.

\setlength{\parindent}{2em}

Let $i$ be the minimal index such that $\nu_{i} \notin \mathcal{P}_{Q}$, noting $1 \leq i \leq n-1$. Let $\omega_{i-1}$ and $\omega_{i}$ be, respectively, the one element of $\pi_{Q}(\nu_{i-1})$ and the one element of $\pi_{Q}(\nu_{i})$, noting that $\omega_{i-1} = \omega_{i}$. According to Theorem \ref{project} there exists a sequence of integers $(n_{j}) \subseteq \{i-1, \ldots, n\}$, containing $i-1$ and $n$, and a corresponding sequence of pants decompositions $\omega_{n_{j}} \in \pi_{Q}(\nu_{n_{j}})$ such that $0 < n_{j+1} - n_{j} \leq 2$, for each $j$, and such that $d_{Q}(\omega_{n_{j}}, \omega_{n_{j+1}}) \leq n_{j+1} - n_{j}$, for each $j$. Necessarily, $\omega_{i-1} = \nu_{i-1}$ and $\omega_{n} = \nu_{n}$. We note that $$d_{Q}(\omega_{i-1}, \omega_{n}) = d_{Q}(\omega_{i}, \omega_{n}) \leq \sum_{j} d_{Q}(\omega_{n_{j}}, \omega_{n_{j+1}}) \leq \sum_{j} n_{j+1} - n_{j} = n - i,$$ and it follows that $$d(\nu_{0}, \nu_{n}) = d(\nu_{0}, \nu_{i-1}) + d(\nu_{i-1}, \nu_{n}) \leq i - 1 + d_{Q}(\omega_{i-1}, \omega_{n}) \leq i - 1 + n - i = n-1.$$ To be more succinct, $d(\nu_{0}, \nu_{n}) \leq n - 1$. This is a contradiction, and the statement of Theorem \ref{main} follows. $\diamond$\\\\

\setlength{\parindent}{0em}

\textbf{$\S5.$ One proof of Corollary \ref{planes}.}\\

\noindent We will need the following two results. The first identifies products of Farey graphs in the pants graph, and the second, stated without proof, recalls a standard property of products of graphs.

\newtheorem{8}[1]{\lemma}

\begin{8}
Let $Q$ be a $2$-handle multicurve on $\Sigma$. Then, the graph $\mathcal{P}_{Q}$ is isomorphic to the product of two Farey graphs.
\end{8}

\noindent \textbf{Proof:} Let $Q_{1}$ and $Q_{2}$ be two disjoint codimension $1$ multicurves on $\Sigma$ such that $Q_{1} \cap Q_{2}$ is equal to $Q$. Then, both $\mathcal{P}_{Q_{1}}$ and $\mathcal{P}_{Q_{2}}$ are subgraphs of $\mathcal{P}_{Q}$ and the product $\mathcal{P}_{Q_{1}} * \mathcal{P}_{Q_{2}}$ is also a subgraph of $\mathcal{P}_{Q}$. Since each vertex $\nu$ of $\mathcal{P}_{Q}$ can be decomposed into the factors $\pi_{Q_{1}}(\nu)$ and $\pi_{Q_{2}}(\nu)$, it follows $\mathcal{P}_{Q}$ is equal to $\mathcal{P}_{Q_{1}} * \mathcal{P}_{Q_{2}}$. Finally, both $\mathcal{P}_{Q_{1}}$ and $\mathcal{P}_{Q_{2}}$ are Farey graphs. $\diamond$\\

\newtheorem{9}[1]{\lemma}

\begin{9}
Let $G_{1}$ and $G_{2}$ be two graphs, and denote by $L_{1}$ and $L_{2}$ convex subgraphs of each respectively. Then, $L_{1} * L_{2}$ is a convex subgraph of $G_{1} * G_{2}$.
\label{someplanes}
\end{9}

\setlength{\parindent}{2em}

A proof of Corollary \ref{planes} can be completed as follows. Let  $Q$ be any $2$-handle multicurve on $\Sigma$, and let $Q_{1}$ and $Q_{2}$ be two disjoint codimension $1$ multicurves such that $Q_{1} \cap Q_{2}$ is equal to $Q$. Let $L_{1}$ be a bi-infinite geodesic line in $\mathcal{P}_{Q_{1}}$, and let $L_{2}$ be a bi-infinite geodesic line in $\mathcal{P}_{Q_{2}}$. Then, $L_{1} * L_{2}$ is a convex subgraph of the totally geodesic subgraph $\mathcal{P}_{Q}$ of $\mathcal{P}(\Sigma)$. Thus, $L_{1} * L_{2}$ is a convex subgraph of $\mathcal{P}(\Sigma)$. Finally, $L_{1} * L_{2}$ is isomorphic to the Cayley graph of $\Bbb{Z} \oplus \Bbb{Z}$ with standard generating set, and as such is a plane. $\diamond$\\~\\

\setlength{\parindent}{0em}

\textbf{References.}\\

\setlength{\parskip}{0.5em plus 0.5em minus 0.5em}

[APS] J. Aramayona, H. Parlier, K. J. Shackleton, \textit{Totally geodesic subgraphs of the pants complex} : arXiv:math.GT/0608752.

[B] J. F. Brock, \textit{The Weil-Petersson metric and volumes of $3$-dimensional hyperbolic convex cores} : Journal of the American Mathematical Society \textbf{16} No. 3 (2003) 495--535.

[BF] J. F. Brock, B. Farb, \textit{Curvature and rank of Teichm\"uller space} : American Journal of Mathematics \textbf{128} (2006) 1--22.

[HT] A. E. Hatcher, W. P. Thurston, \textit{A presentation for the mapping class group of a closed orientable surface} : Topology \textbf{19} (1980) 221--237.

[Mar] D. Margalit, \textit{Automorphisms of the pants complex} : Duke Mathematical Journal \textbf{121} No. 3 (2004) 457--479.

[Mas] H. A. Masur, \textit{Extension of the Weil-Petersson metric to the boundary of Teichmuller space} : Duke Mathematical Journal \textbf{43} no. 3 (1976) 623--635.

[MasMi] H. A. Masur, Y. N. Minsky, \textit{Geometry of the complex of curves II: Hierarchical structure} : Geometry \& Functional Analysis \textbf{10} (2000) 902--974.

[MasS] H. A. Masur, S. Schleimer, \textit{The pants complex has only one end} : in ``Spaces of Kleinian groups'' (eds. Y. N. Minsky, M. Sakuma, C. M. Series) London Math. Soc. Lecture Note Ser. \textbf{329} (2006) 209--218.

[W] S. A. Wolpert, \textit{Geometry of the Weil-Petersson completion of Teichm\"uller space} : Surveys in Differential Geometry VIII: Papers in honor of Calabi, Lawson, Siu and Uhlenbeck, editor S. T. Yau, International Press (2003).\\

\newpage

\setlength{\parskip}{0em plus 0.5em minus 0.5em}

\begin{small}
Javier Aramayona\\
CMI\\
Universit\'e de Provence\\
39, rue Joliot Curie\\
13453 Marseille\\
France\\
homepage: http://www.maths.warwick.ac.uk/$\sim$jaram\\
e-mail: jaram@maths.warwick.ac.uk\\

Hugo Parlier\\
Section de Math\'ematiques\\
Universit\'e de Gen\`eve\\
1211 Gen\`eve 4\\
Suisse\\
homepage: http://www.unige.ch/math/folks/parlier/\\
e-mail: hugo.parlier@math.unige.ch\\

Kenneth J. Shackleton (corresponding author)\\
(Professor Sadayoshi Kojima Laboratory)\\
Department of Mathematical and Computing Sciences\\
Tokyo Institute of Technology\\
2-12-1 O-okayama\\
Meguro-ku\\
Tokyo 152-8552\\
Japan\\
homepage: http://www.maths.soton.ac.uk/$\sim$kjs\\
e-mail: shackleton.k.aa@m.titech.ac.jp\\
e-mail: kjs2006@alumni.soton.ac.uk
\end{small}

\end{document}